# Existence of Solutions to Split Variational Inequality Problems and Split Minimization Problems in Banach Spaces

Jinlu Li

Department of Mathematical Sciences
Shawnee State University
Portsmouth, Ohio 45662
U.S.A.

**Abstract:** In this paper, we extend the concept of split variational inequality problems from Hilbert spaces to Banach spaces. Then we apply the Fan-KKM theorem to prove the existence of solutions to some split variational inequality problems and some split minimization problems in Banach spaces. We also apply fixed point theorems on Chain complete posets to show the solvability of some split variational inequality problems in partially ordered Banach space.

**Keywords:** Fan-KKM theorem; split variational inequality problem in Banach space; split minimization problem; fixed point; partially ordered Banach space

**2010 Mathematics Subject Classification:** 49J40; 49J52; 49J20

## 1. Introduction

Split variational inequality problems in Hilbert spaces have been studied by many authors (see [3], [6-9] and [14]). This topic has been caused many authors attention. It is because that, excepting the directed applications to some medical image reconstruction problems, split variational inequality problems are natural extensions of the classical variational inequality problems. It is also very closely connected with the theory of split feasibility problems (see [3-6]). Recently, research on split variational inequality problems has become a new trend of study in nonlinear analysis.

Concurrent to the development of split variational inequality theory, some researchers are studying the theory of split variational inclusion problems and split feasibility problems in Hilbert spaces and Banach spaces (for split feasibility problems in Hilbert spaces, see [1-2], [4-5], [14-15]; for split feasibility problems in Banach spaces, see [17]). The theory of split feasibility problems stems from some medical image reconstruction problems, which was transformed into some modeling inverse problems in finite-dimensional Hilbert spaces by Censor and Elfving in 1994 (see [4], [18-19]). Since then, split feasibility problems together with split variational inequality problems have attracted more attention.

We recall the concept of split variational inequality problems in Hilbert spaces below.

Let $H_1$ and $H_2$ be Hilbert spaces. For any given nonempty, closed and convex subsets $C \subseteq H_1$ and $D \subseteq H_2$, given operators $f: C \to H_1$, $g: D \to H_2$ and a bounded linear operator $A: H_1 \to H_2$, the split variational inequality problem associated with $f, g, A, C, D$ is formulated as follows:

$$\text{find } x_* \in C \text{ such that } \langle f(x_*), x - x_* \rangle \geq 0 \text{ for all } x \in C, \tag{1.1}$$

$$\text{and such that the point } y_* = Ax_* \in D \text{ solves } \langle g(y_*), y - y_* \rangle \geq 0 \text{ for all } y \in D. \tag{1.2}$$

In the literature of the research on split variational inequality problems in Hilbert spaces, such as [3], [6-9] and [14], the authors set certain conditions for the operators $f$, $g$ and $A$, such that they developed some algorithms to create iterative schemes that converge, weakly or strongly, to the solutions of the considered split variational inequality problems defined in (1.1) and (1.2). The essential prerequisite to constructing such convergent sequences is to assume the existence of solutions to these problems. An important question is immediately raised: what are the conditions on $f$, $g$ and $A$ to assure the existence of solutions to the split variational inequality problems defined in (1.1) and (1.2)? The goal of this paper is to try to answer this question. More precisely, in this paper, we will study the existence of solutions to split variational inequality problems in Banach spaces. It is clearly to see that the linear operator $A$ plays very important role in these problems. So we investigate some applicable conditions about $A$, together with $f$ and $g$ to insure the solvability of these problems in Banach spaces and in partially ordered Banach spaces.

This paper is organized as follows: in section 2, we recall the concept of KKM mappings and the Fan-KKM theorem and prove a solution existence theorem for some split variational inequality problems in Banach spaces; in section 3, we recall some concepts of partially ordered vector spaces and apply some fixed point theorems on chain complete poset to show the solvability of some split variational inequality problem in partially ordered Banach spaces; In section 4, we prove a theorem for the existence of solutions to some split convex minimization problems.

## 2. Split variational inequality problems in Banach spaces

### 2.1. Definitions of split variational inequality problems in Banach spaces

Let $X$ and $Y$ be Banach spaces with dual spaces $X^*$ and $Y^*$, respectively. Without any confusion caused, we use $\langle \cdot, \cdot \rangle$ to denote the pairing between $X^*$ and $X$ and between $Y^*$ and $Y$. For any given nonempty, closed and convex subsets $C \subseteq X$ and $D \subseteq Y$, given operators $f: C \to X^*$, $g: D \to Y^*$ and a bounded linear operator $A: X \to Y$, the split variational inequality problem associated with $f, g, A, C, D$ (denoted by SVIP($f, C, A, g, D$)) is formulated as follows:

$$\text{find } x_* \in C \text{ such that } \langle f(x_*), x - x_* \rangle \geq 0 \text{ for all } x \in C, \tag{2.1}$$

$$\text{and such that the point } y_* = Ax_* \in D \text{ solves } \langle g(y_*), y - y_* \rangle \geq 0 \text{ for all } y \in D. \tag{2.2}$$

Such a pair $(x_*, y_*) = (x_*, Ax_*)$ is called a solution to SVIP($f, C, A, g, D$). Let $\mathcal{S}(f, C, A, g, D)$ denote the set of all solutions to SVIP($f, C, A, g, D$).

When looked at separately, the problem (2.1) is the classic variational inequality problem (VIP) on Banach spaces, which is usually denoted by VIP($f$, $C$). On the other hand, when we take a special case in SVIP($f$, $C$, $A$, $g$, $D$), such as: $X = Y$, $C = D$, $f = g$ and $A = I$, that is the identity mapping on $X$, then the split variational inequality problem SVIP($f$, $C$, $A$, $g$, $D$) becomes the classic variational inequality problem VIP($f$, $C$). So we can consider split variational inequality problems to be the natural extensions of classic variational inequality problems in Banach spaces.

## 2.2 Convexity direction reserving mappings and examples

The existence of solutions to the classic variational inequality problems VIP($f$, $C$) and VIP($g$, $D$) in Banach spaces has been studied by many authors, where the mappings $f$ and $g$ are required to satisfy a certain type of conditions. For split variational inequality problems, it is obvious that the linear operator $A$ plays a very important role in the existence of solutions to SVIP($f$, $C$, $A$, $g$, $D$). One of the goals in this paper is to investigate some applicable conditions for $A$, $f$, and $g$ to assure the existence of solutions to SVIP($f$, $C$, $A$, $g$, $D$). For this purpose, we introduce some concepts below.

For any given vectors $x_1, x_2, \cdots, x_n$ of a linear vector space with an arbitrary positive integer $n$, and $0 < \lambda_1, \lambda_2, \cdots, \lambda_n < 1$ satisfying $\Sigma_{1 \leq i \leq n} \lambda_i = 1$, the vector $\Sigma_{1 \leq i \leq n} \lambda_i x_i$ in the same vector space is called a strictly linear combination of the vectors $x_1, x_2, \cdots, x_n$.

Let $X$, $Y$, $C$, $D$, $f$ and $g$ be given as above. A bounded linear operator $A: X \to Y$ is said to have the *convexity direction reserving property* (or it is said to be *convexity direction reserved*) with respect to the mappings $f$ and $g$ on $C$ and $D$ if, for any given vectors $x_1, x_2, \cdots, x_n$ in $C$, and for any strictly linear combination $x' = \Sigma_{1 \leq i \leq n} \lambda_i x_i$, there is a positive integer $j \leq n$ such that the following two inequalities simultaneously hold:

$$\langle f(x'), x_j - x' \rangle \geq 0 \text{ and } \langle g(Ax'), Ax_j - Ax' \rangle \geq 0 \tag{2.3}$$

To understand the meaning of the convexity direction reserving property of a bounded linear operator $A$, we consider the case that $n = 2$. For arbitrary two points $x_1, x_2$ in $C$, and for any strictly linear combination $x' = \lambda x_1 + (1 - \lambda) x_2$, for some $0 < \lambda < 1$, since $x_1 - x'$ and $x_2 - x'$ have the opposite directions, then we must have

$$\langle f(x'), x_1 - x' \rangle \geq 0 \text{ or } \langle f(x'), x_2 - x' \rangle \geq 0.$$

Similarly, since $Ax_1 - Ax'$ and $Ax_2 - Ax'$ have the opposite directions, we have

$$\langle g(Ax'), Ax_1 - Ax' \rangle \geq 0 \text{ or } \langle g(Ax'), Ax_2 - Ax' \rangle \geq 0.$$

The convexity direction reserving property of $A$ demonstrates that there is a number $j = 1$ or $2$ such that the two inequalities in (2.3) simultaneously hold. We provide some examples for direction reserved mappings below.

**Example 2.1**. Let $X = Y$, $C = D$, where $C$ contains the origin point, and let $f = g$ be linear operators. For any given positive number $\alpha \leq 1$, define $Ax = \alpha x$, for all $x \in X$. Then $A$ is a bounded linear operator on $X$ that has the convexity direction reserving property with respect to the mapping $f$ on $C$.

In particular, if $C$ is a closed convex cone in $X$, then for any $\beta > 0$, the operator $A$ defined by $Ax = \beta x$, has the convexity direction reserving property with respect to the mapping $f$ on $C$.

Since many useful Banach spaces with their dual spaces have biorthogonal systems (see [16]), next we consider bounded linear operators on such Banach spaces.

**Example 2.2**. Let $X$ and $Y$ be two Banach spaces with dual spaces $X^*$ and $Y^*$, respectively. Suppose that $X$ has a basis $\{d_k\}$ and $X^*$ has a basis $\{u_k\}$ such that $(d_k, u_k)$ is a biorthogonal system (see [16]); that is,

$$\langle u_k, d_l \rangle = \delta_{kl}, \text{ for } k, l = 1, 2, \cdots.$$

Let $f: C \to X^*$ be a continuous operator defined by

$$f(x) = \Sigma_{1 \leq k < \infty} f_k(x) u_k, \text{ for all } x \in C,$$

where $\{f_k\}$ is a sequence of continuous functions defined on $C$, such that, for any $s \in X$ with $s = \Sigma_{1 \leq k < \infty} b_k d_k$, we have

$$\langle f(x), s \rangle = \Sigma_{1 \leq k < \infty} b_k f_k(x).$$

Let $x_1, x_2, \cdots, x_n$ be arbitrarily given vectors in $C$, for some positive integer $n$, with

$$x_i = \Sigma_{1 \leq k < \infty} b_{ik} d_k, i = 1, 2, \cdots, n.$$

Let $x' = \Sigma_{1 \leq i \leq n} \lambda_i x_i$ be a strictly linear combination of arbitrarily given vectors $x_1, x_2, \cdots, x_n$. Then, for any fixed $j = 1, 2, \cdots, n$, we have

$$\begin{aligned} x_j - x' &= \Sigma_{i \neq j} \lambda_i (x_j - x_i) \\ &= \Sigma_{i \neq j} \lambda_i \Sigma_{1 \leq k < \infty} (b_{jk} - b_{ik}) d_k \\ &= \Sigma_{1 \leq k < \infty} (\Sigma_{i \neq j} \lambda_i (b_{jk} - b_{ik})) d_k. \end{aligned}$$

It follows that

$$\langle f(x'), x_j - x' \rangle = \Sigma_{1 \leq k < \infty} (\Sigma_{i \neq j} \lambda_i (b_{jk} - b_{ik})) f_k(x'). \tag{2.4}$$

Suppose that $Y$ and $Y^*$ have bases $\{e_k\}$ and $\{v_k\}$, respectively, such that $(e_k, v_k)$ is a biorthogonal system (see Theorems 17.1 and 17.2 in [16]). Let $A: X \to Y$ be a bounded linear operator such that

$$Ax = \Sigma_{1 \leq k < \infty} \langle u_k, x \rangle e_k, \text{ for all } x \in X.$$

Let $g: D \to Y^*$ be a continuous operator defined by

$$g(y) = \Sigma_{1 \leq k < \infty} g_k(y) v_k, \text{ for all } y \in C,$$

where $\{g_k\}$ is a sequence of continuous functions defined on $D$, such that, for any $t \in Y$ with $t = \Sigma_{1 \leq k < \infty} c_k e_k$, we have

$$\langle g(y), t \rangle = \Sigma_{1 \leq k < \infty} c_k g_k(y).$$

For the bounded linear operator $A: X \to Y$ defined as above, we calculate

$$\begin{aligned}
&Ax_j - Ax' \\
&= \Sigma_{i \neq j} \lambda_i A(x_j - x_i) \\
&= \Sigma_{i \neq j} \lambda_i \Sigma_{1 \leq k < \infty} \langle u_k, x_j - x_i \rangle e_k \\
&= \Sigma_{1 \leq k < \infty} (\Sigma_{i \neq j} \lambda_i \langle u_k, x_j - x_i \rangle) e_k.
\end{aligned}$$

Then we obtain

$$\begin{aligned}
&\langle g(Ax'), Ax_j - Ax' \rangle \\
&= \Sigma_{1 \leq k < \infty} (\Sigma_{i \neq j} \lambda_i \langle u_k, x_j - x_i \rangle) g_k(Ax') \\
&= \Sigma_{1 \leq k < \infty} (\Sigma_{i \neq j} \lambda_i \langle u_k, \Sigma_{1 \leq l < \infty} (b_{jl} - b_{il}) d_l \rangle) g_k(Ax') \\
&= \Sigma_{1 \leq k < \infty} (\Sigma_{i \neq j} \lambda_i (b_{jk} - b_{ik})) g_k(Ax').
\end{aligned} \quad (2.5)$$

Suppose that $f: C \to X^*$ and $g: D \to Y^*$ are "positive linearly dependent" with respect to this bounded linear operator $A: X \to Y$; that is, there is a positive number $\alpha$ such that

$$g_k(Ax) = \alpha f_k(x), \text{ for } k = 1, 2, \cdots . \quad (2.6)$$

Then from (2.6) and by combining (2.4) and (2.5), $A$ has the convexity direction reserving property with respect to the mappings $f$ and $g$ on $C$ and $D$.

### 2.3. Fan-KKM Theorem

Fan-KKM theorem is applied in the proof of the main theorem in this section. We recall the definition of KKM mappings and Fan-KKM theorem below. For more details regarding to KKM mappings and Fan-KKM theorem, please refer to Fan [10] and Park [15].

Let $K$ be a nonempty subset of a linear space $B$. A set-valued mapping $T: K \to 2^B \setminus \{\emptyset\}$ is said to be a KKM mapping if for any finite subset $\{y_1, y_2, \ldots, y_n\}$ of $K$, we have

$$\text{co}\{y_1, y_2, \ldots, y_n\} \subseteq \cup_{1 \leq i \leq n} T(y_i),$$

where $\text{co}\{y_1, y_2, \ldots, y_n\}$ denotes the convex hull of $\{y_1, y_2, \ldots, y_n\}$.

**Fan-KKM Theorem**. *Let K be a nonempty convex subset of a Hausdorff topological vector space B and let T: $K \to 2^B \setminus \{\emptyset\}$ be a KKM mapping with closed values. If there exists a point $y_0 \in K$ such that $T(y_0)$ is a compact subset, then*

$$\bigcap_{y \in K} T(y) \neq \emptyset.$$

## 2.4 Existence of solutions to split variational inequality problems in Banach spaces

Now we state and prove the main theorem in this section.

**Theorem 2.3**. *Let X and Y be Banach spaces and $C \subseteq X$ and $D \subseteq Y$ be nonempty, closed and convex subsets. Let $f: C \to X^*$, $g: D \to Y^*$ be continuous operators and let $A: X \to Y$ be a bounded linear operator satisfying the following conditions*:

(a1) $AC = D$;
(a2) *A is convexity direction reserved with respect to f and g on C and D*.

*If there is $(x_0, Ax_0) \in C \times D$ such that*

$$\{(x', Ax') \in C \times D: \langle f(x'), x_0 - x' \rangle \geq 0 \text{ and } \langle g(Ax'), Ax_0 - Ax' \rangle \geq 0\} \text{ is compact,} \quad (2.7)$$

*then SVIP(f, C, A, g, D) has a solution*.

**Proof**. Denote the graph of the operator $A$ by

$$M = \{(x, Ax) \in C \times D: x \in C\}.$$

Since $A$ is a bounded linear operator from the Banach space $X$ to Banach space $Y$, the convexity of $C$ and $D$ implies that $M$ is convex. By the closed graph theorem in Banach spaces, from the continuity (bounded) of the linear operator $A$, it yields that the graph $M$ of $A$ is closed with the product topology. Hence $M$ is a nonempty, closed and convex subset of $C \times D \subseteq X \times Y$. We define the product mapping of the mappings $f$ and $g$ by $h: C \times D \to X^* \times Y^*$ as

$$h(x, y) = (f(x), g(y)), \text{ for every } (x, y) \in C \times D,$$

such that for every $(x, y) \in C \times D$ and for every $(s, t) \in X \times Y$,

$$h(x, y)(s, t) = (f(x), g(y))(s, t) = (\langle f(x), s \rangle, \langle g(y), t \rangle) \in R_2. \quad (2.8)$$

Since $f(x) \in X^*$ and $g(y) \in Y^*$, then $h$ is an operator from $C \times D$ to $\mathcal{L}(X \times Y, R^2)$. Then, by the mapping $h$ in (2.8), the split variational inequality problem SVIP(f, C, A, g, D) by (2.1) and (2.2) can be transformed into the following ordered variational inequality problem OVIP(h, C×D):

find a point $(x_*, Ax_*) \in C \times D$ such that $h(x_*, Ax_*)((x, y) - (x_*, Ax_*)) \succcurlyeq_2 0$, for all $(x, y) \in C \times D$,

where the component partial order $\succcurlyeq_2$ on $R_2$ is defined by: for any $(t_1, s_1)$, $(t_2, s_2) \in R_2$,

$$(t_2, s_2) \succcurlyeq_2 (t_1, s_1) \text{ if and only if } t_2 \geq t_1 \text{ and } s_2 \geq s_1.$$

It is well-known that the partial order $\succcurlyeq_2$ on $R_2$ is a lattice and $(R_2, \succcurlyeq_2)$ is a 2-d Hilbert lattice. To prove this theorem by applying the Fan-KKM Theorem, define a mapping $T: M \to 2^M$ by

$$T(x, Ax) = \{(x', Ax') \in M: h(x', Ax')((x, Ax) - (x', Ax')) \succcurlyeq_2 0\}, \text{ for all } (x, Ax) \in M.$$

Since $(x, Ax) \in T(x, Ax)$, from the continuity of the mappings $A$, $f$ and $g$, it yields that, for every $(x, Ax) \in M$, $T(x, Ax)$ is a nonempty and closed subset of $M \subseteq C \times D \subseteq X \times Y$.

Next we show that $T$ is a KKM mapping. For an arbitrary positive integer $n$ greater than 1, we arbitrarily take points $(x_1, Ax_1)$, $(x_2, Ax_2)$, $\cdots$, $(x_n, Ax_n) \in M$ and numbers $0 < \lambda_1, \lambda_2, \cdots, \lambda_n < 1$ such that $\Sigma_{1 \leq i \leq n} \lambda_i = 1$. Let

$$(x', Ax') = \Sigma_{1 \leq i \leq n} \lambda_i (x_i, Ax_i), \text{ where } x' = \Sigma_{1 \leq i \leq n} \lambda_i x_i \text{ and } Ax' = \Sigma_{1 \leq i \leq n} \lambda_i Ax_i.$$

Then $x' = \Sigma_{1 \leq i \leq n} \lambda_i x_i$ is a strictly linear combination of the vectors $x_1, x_2, \cdots, x_n$. From condition a2, $A$ has the convexity direction reserving property with respect to $f$ and $g$ on $C$ and $D$. It follows that there is a positive integer $j \leq n$ such that

$$\langle f(x'), x_j - x' \rangle \geq 0 \text{ and } \langle g(Ax'), Ax_j - Ax' \rangle \geq 0. \tag{2.3}$$

That is,

$$h(x', Ax')((x_j, Ax_j) - (x', Ax')) \succcurlyeq_2 0.$$

It implies $(x', Ax') \in T(x_j, Ax_j)$, for some $j = 1, 2, \cdots, n$. We obtain

$$(x', Ax') = \Sigma_{1 \leq i \leq n} \lambda_i (x_i, Ax_i) \in T(x_j, Ax_j) \subseteq \cup_{1 \leq i \leq n} T(x_i, Ax_i).$$

Hence, $T$ is a KKM mapping. For the point $(x_0, Ax_0) \in M$ given in condition (2.7) in this theorem, we see that $T(x_0, Ax_0)$ is compact. By applying Fan-KKM theorem, $\cap_{(x, Ax) \in M} T(x, Ax) \neq \varnothing$. Then taking any $(x_*, Ax_*) \in \cap_{(x, Ax) \in M} T(x, Ax)$, it satisfies

$$h(x_*, Ax_*)((x, Ax) - (x_*, Ax_*)) \succcurlyeq_2 0, \text{ for all } (x, Ax) \in M. \tag{2.9}$$

The order-inequality (2.9) is equivalent to

$$\langle f(x_*), x - x_* \rangle \geq 0 \text{ for all } x \in C, \tag{2.10}$$
and

$$\langle g(Ax_*), Ax - Ax_* \rangle \geq 0 \text{ for all } x \in C. \tag{2.11}$$

Let $y_* = Ax_*$. From condition (a1) of the mapping $A$ (notice that $x_*$ and $y_* = Ax_*$ in (2.10) and (2.11) are fixed), (2.11) reduces to

$$\langle g(y_*), y - y_* \rangle \geq 0 \text{ for all } y \in D. \tag{2.12}$$

Combining (2.10), (2.12) and $y_* = Ax_*$, it implies that $(x_*, y_*) = (x_*, Ax_*)$ is a solution of the split variational inequality problem SVIP($f, C, A, g, D$). □

## 3. Split variational inequality problems in partially ordered Banach spaces

In this section, we consider split variational inequality problems in Banach spaces equipped with some partial orders. The structures of the partial orders with their properties on the considered Banach spaces will help us not only to prove the existence of solutions to some split variational inequality problems, but also to study the inductive properties of the solution sets.

### 3.1. Some preliminaries of posets and partially topological vector spaces

When we study fixed point theory on posets, the order-monotonic property of the considered mapping and the chain complete property of the underlying spaces are important for application. We recall some concepts and properties of posets and partially topological vector spaces below.

Let $(U, \succcurlyeq^U)$ and $(V, \succcurlyeq^V)$ be posets. A set-valued mapping $T: U \to 2^V \setminus \{\emptyset\}$ is said to be isotone or order-increasing upward whenever $u_1 \preccurlyeq^U u_2$ in $U$ implies that, for any $z \in T(u_1)$, there is a $w \in T(u_2)$ such that $z \preccurlyeq^V w$.

Let $X$ be a topological vector space equipped with a partial order $\succcurlyeq^X$. The topology on $X$ is said to be natural with respect to the given partial order $\succcurlyeq^X$ if, for any $u \in X$, the following $\succcurlyeq^X$-intervals are closed:

$$(u] = \{x \in X : x \preccurlyeq^X u\} \text{ and } [u) = \{x \in X : x \succcurlyeq^X u\}.$$

If the topology on $X$ is natural with respect to $\succcurlyeq^X$ and for any $u, v \in X$ satisfying $v \succcurlyeq^X u$, the following $\succcurlyeq^X$-inequalities hold:

$$\alpha v \succcurlyeq^X \alpha u \text{ and } v + w \succcurlyeq^X u + w, \text{ for any } \alpha \geq 0 \text{ and for any } w \in X,$$

then $X$ is called a partially ordered vector space with $\succcurlyeq^X$ and it is written as $(X, \succcurlyeq^X)$.

**Observation 3.1.** If $X$ is a Banach space equipped with a partial order $\succcurlyeq^X$, then the norm topology is natural with respect to the partial order $\succcurlyeq^X$, if and only if, the weak topology on $X$ is natural with respect to the partial order $\succcurlyeq^X$. Hence $(X, \succcurlyeq^X)$ is a partially ordered Banach space with the strong topology, if and only if, it is with the weak topology.

The chain complete property of the underlying space in fixed point theory on posets play an important role for the existence of fixed point of a considered mapping. We provide some useful examples of chain complete posets below, which is a consequence of Lemma 5.2 in [13].

**Lemma 3.2**.
  (i) Let $C$ be a nonempty compact subset of a partially ordered Banach space $(X, \succcurlyeq^X)$. Then $(C, \succcurlyeq^X)$ is $\succcurlyeq^X$-chain complete.
  (ii) Let $C$ be a nonempty bounded, closed and convex subset of a partially ordered reflexive Banach space $(X, \succcurlyeq^X)$. Then $(C, \succcurlyeq^X)$ is $\succcurlyeq^X$-chain complete.

**Proof**. Part (i) immediately follows from Lemma 5.2 in [13]. Since every nonempty bounded, closed and convex subset of a reflexive Banach space is weakly compact, then, from the Observation 3.1 and Lemma 5.2 in [13], Part (ii) is obtained. □

The concept of universally inductive posets is introduced in [12] and [13], which has been used in the proof of fixed point existence theorems. We recall it here. A nonempty subset $A$ of a poset $(P, \succcurlyeq)$ is said to be universally inductive in $P$ whenever, any given chain $\{x_\alpha\} \subseteq P$ satisfying that if every element $x_\beta \in \{x_\alpha\}$ has an upper cover in $A$, then the chain $\{x_\alpha\}$ has an upper bound in $A$. We find that the definition of universally inductive posets is a very broad concept. It includes many useful subsets in posets and in partially ordered topological spaces. We list some of them below, which are from [12] and [13].

**Lemma 3.2 [12]**. *Every inductive subset A in a chain complete poset with a finite number of maximal elements is universally inductive.*

**Lemma 5.3 [13]**. *Every non-empty compact subset of a partially ordered Hausdorff topological space is universally inductive.*

**Lemma 5.4 [13]**. *Every non-empty bounded, closed and convex subset of a partially ordered reflexive Banach space is universally inductive.*

The proof of the main theorem for the existence of solutions to some split variational inequality problems in partially ordered Banach spaces is based on the following theorem that is from [13].

**Theorem 3.1 [13]**. *Let $(P, \succcurlyeq)$ be a chain complete poset and let $T: P \to 2^P \setminus \{\emptyset\}$ be a set-valued mapping satisfying the following three conditions*:

  A1. *T is order-increasing upward*;
  A2. *$(T(x), \succcurlyeq)$ is universally inductive, for every $x \in P$*;
  A3. *There is an element $x_0$ in $P$ and $v_0 \in T(x_0)$ with $x_0 \preccurlyeq v_0$.*
*Then*
  (i)  *$(\mathcal{F}(T), \succcurlyeq)$ is a nonempty inductive poset*;
  (ii) *$(\mathcal{F}(T) \cap [x_0), \succcurlyeq)$ is a nonempty inductive poset; and T has an $\succcurlyeq$-maximal fixed point $x^*$ with $x^* \succcurlyeq x_0$.*

Where $\mathcal{F}(T)$ denotes the set of fixed points of $T$.

## 3.2. Minimizing acceptances of operators and their properties

Let $C \subseteq X$ and $D \subseteq Y$ be nonempty subsets of Banach spaces, $X$ and $Y$, respectively. Given operators $f: C \to X^*$, $g: D \to Y^*$, we define set-valued mappings $F: C \to 2^C$ and $G: D \to 2^D$ as follows:

$$F(x) = \{u \in C: \langle f(x), u \rangle = \min_{t \in C} \langle f(x), t \rangle\}, \text{ for all } x \in C,$$

and
$$G(y) = \{v \in D: \langle g(y), v \rangle = \min_{s \in C} \langle g(y), s \rangle\}, \text{ for all } y \in D.$$

$F$ and $G$ are called the *minimizing acceptances* of $f$ and $g$, respectively.

We list some properties of $F$ and $G$ below, which will be useful in the sequel of this paper.

**Lemma 3.3**. Let $(X, \succcurlyeq^X)$ and $(Y, \succcurlyeq^Y)$ be partially ordered Banach spaces and let $C \subseteq X$ and $D \subseteq Y$ be nonempty compact subsets. Let $f: C \to X^*$, $g: D \to Y^*$ be continuous operators. Then the values of their minimizing acceptances $F(x)$, for every $x \in C$, $G(y)$, for every $y \in D$, are nonempty and compact subsets in $X$ and $Y$, respectively. Furthermore, both $(F(x), \succcurlyeq^X)$ and $(G(y), \succcurlyeq^Y)$ are chain complete.

**Proof**. The proof of the first part is straightforward. Then by the first part and Lemma 3.2, the second part follows immediately. □

**Definition 3.4**. Let $(X, \succcurlyeq^X)$ be a partially ordered Banach space and let $C \subseteq X$ be nonempty subset. An operator $f: C \to X^*$ is said to be $\succcurlyeq^X$-*decrement* if, for any $x_1 \preccurlyeq^X x_2$ in $C$, and for any points $u_1, w \in C$ satisfy

$$\langle f(x_1), u_1 \rangle \leq \langle f(x_1), w \rangle,$$

then there is $u_2 \in C$ such that $u_1 \preccurlyeq^X u_2$ and

$$\langle f(x_2), u_2 \rangle \leq \langle f(x_2), w \rangle.$$

**Lemma 3.5**. Let $(X, \succcurlyeq^X)$ be a partially ordered Banach space and let $C \subseteq X$ be nonempty subset. Let $f: C \to X^*$ be an $\succcurlyeq^X$-decrement operator. Then its minimizing acceptance $F$ is $\succcurlyeq^X$-increasing upward.

**Proof**. For any given $x_1 \preccurlyeq^X x_2$ in $C$ and for any $u_1 \in F(x_1)$, we take an arbitrary point $w \in F(x_2)$, which satisfies

$$\langle f(x_1), u_1 \rangle = \min_{t \in C} \langle f(x), t \rangle \leq \langle f(x_1), w \rangle.$$

Since the operator $f: C \to X^*$ is $\succcurlyeq^X$-decrement, there is $u_2 \in C$ such that $u_1 \preccurlyeq^X u_2$ and

$$\langle f(x_2), u_2 \rangle \leq \langle f(x_2), w \rangle = \min_{t \in C} \langle f(x), t \rangle.$$

It implies $u_2 \in F(x_2)$. □

**Lemma 3.6.** Let $(X, \succcurlyeq^X)$ be a partially ordered Banach space and let $C \subseteq X$ be nonempty subset. Let $F$ be the minimizing acceptance of an $\succcurlyeq^X$-decrement operator $f: C \to X^*$. Then, for every $u \in F(x)$, with any given $x \in C$, we have

$$[u) = \{w \in C: w \succcurlyeq^X u\} \subseteq F(x).$$

**Proof.** For any given $x \in C$ and for an arbitrary point $u \in F(x)$, take any $v \in C$ with $v \succcurlyeq^X u$. Then, for any $w \in C$, we have

$$\langle f(x), u \rangle \leq \langle f(x), w \rangle.$$

Since $x \succcurlyeq^X x$ and $v \succcurlyeq^X u$, from the $\succcurlyeq^X$-decrement property of $f$, it implies

$$\langle f(x), v \rangle \leq \langle f(x), w \rangle, \text{ for any } w \in C.$$

It follows that $v \in F(x)$. □

**Definition 3.7.** A bounded linear operator $A: X \to Y$ is said to be order-pseudomonotone with respect to operators $f: C \to X^*$, $g: D \to Y^*$ if, for any given $x$ in $C$, and for any points $u, w \in C$ with $u \not\prec^X w$ ($w \not\prec^X u$), then,

$$\langle f(x), u \rangle \leq \langle f(x), w \rangle \text{ implies } \langle g(Ax), Au \rangle \leq \langle g(Ax), Aw \rangle. \tag{3.1}$$

### 3.3. Existence of solutions to split variational inequality problems in partially ordered Banach spaces

Let $(X, \succcurlyeq^X)$ and $(Y, \succcurlyeq^Y)$ be partially ordered Banach spaces. Let $\succcurlyeq^{X \times Y}$ denote the component partial order on the product space of $(X, \succcurlyeq^X)$ and $(Y, \succcurlyeq^Y)$ that is defined as follows: for any points $(x_1, y_1)$, $(x_2, y_2) \in X \times Y$, we say that

$$(x_2, y_2) \succcurlyeq^{X \times Y} (x_1, y_1) \text{ if and only if } x_2 \succcurlyeq^X x_1 \text{ and } y_2 \succcurlyeq^Y y_1.$$

It follows that $(X \times Y, \succcurlyeq^{X \times Y})$ is also a partially ordered Banach space.

**Theorem 3.8.** *Let $(X, \succcurlyeq^X)$ and $(Y, \succcurlyeq^Y)$ be partially ordered Banach spaces and let $C \subseteq X$ and $D \subseteq Y$ be nonempty compact subsets. Let $f: C \to X^*$, $g: D \to Y^*$ be continuous operators and let $A: X \to Y$ be a bounded linear operator satisfying the following conditions*:

(a0) $f$ is $\geqslant^X$-decrement and $g$ is $\geqslant^Y$-decrement;
(a1) $AC = D$;
(b1) $A$ is order-increasing;
(b2) $A$ is order-pseudomonotone with respect to $f$ and $g$.

*Suppose that there is $(x_0, Ax_0) \in C \times D$ such that*

$$\text{there exists } (w, Aw) \in F(x_0) \times G(Ax_0) \text{ satisfying } x_0 \leqslant^X w. \tag{3.2}$$

*Then* SVIP$(f, C, A, g, D)$ *has a solution. Moreover, we have*

(i) $(\mathcal{S}(f, C, A, g, D), \geqslant^{X \times Y})$ *is a nonempty and inductive poset*;
(ii) $(\mathcal{S}(f, C, A, g, D) \cap [(x_0, Ax_0)), \geqslant^{X \times Y})$ *is also a nonempty and inductive poset*.

**Proof**. Similarly to the proof of Theorem 2.3, we write the graph of the operator $A$ by

$$M = \{(x, Ax) \in C \times D : x \in C\}.$$

Since $A$ is linear and continuous, then $M$ is a closed and convex subset in $C \times D \subseteq X \times Y$. Hence $M$ is a compact subset in $C \times D$ with respect to the product topology. From Lemma 3.2, $(M, \geqslant^{X \times Y})$ is an $\geqslant^{X \times Y}$-chain complete poset. On $(M, \geqslant^{X \times Y})$, a set-valued mapping $T: M \to 2^M$ is defined as

$$T(x, Ax) = \{(x', Ax') \in M : x' \in F(x) \text{ and } Ax' \in G(Ax)\}, \text{ for every } (x, Ax) \in M. \tag{3.3}$$

Since $C$ and $D$ are compact and $f$ and $g$ are continuous, from Lemma 3.3, the minimizing acceptances $F(x)$ and $G(y)$ are nonempty and compact subsets in $X$ and $Y$, for every $x \in C$ and $y \in D$, respectively. Furthermore, both $(F(x), \geqslant^X)$ and $(G(y), \geqslant^Y)$ are chain complete. It implies that $(T(x, Ax), \geqslant^{X \times Y})$ is an $\geqslant^{X \times Y}$-chain complete subset in $(M, \geqslant^{X \times Y})$, for all $(x, Ax) \in M$.

Next we show that $T(x, Ax) \neq \emptyset$, for every $(x, Ax) \in M$. Since $(F(x), \geqslant^X)$ is $\geqslant^X$-chain complete, so it is inductive. Then there is an $\geqslant^X$-maximal element $u$ in $(F(x), \geqslant^X)$ such that

$$\langle f(x), u \rangle \leq \langle f(x), w \rangle, \text{ for all } w \in C. \tag{3.4}$$

Since $C$ is a nonempty compact subset in Banach space $X$, then from Lemma 3.2, $(C, \geqslant^X)$ is also $\geqslant^X$-chain complete; and therefore, it is $\geqslant^X$-inductive. Hence $(C, \geqslant^X)$ has $\geqslant^X$-maximal elements. From Lemma 3.6, it follows that any $\geqslant^X$-maximal element in $(F(x), \geqslant^X)$ must be an $\geqslant^X$-maximal element in $(C, \geqslant^X)$. It implies that $w \not\succ^X u$, for any point $w \in C$. Since the operator $A: X \to Y$ is order-pseudomonotone with respect to operators $f$ and $g$, then from (3.4) and applying (3.1), the following order-inequality holds:

$$\langle g(Ax), Au \rangle \leq \langle g(Ax), Aw \rangle, \text{ for all } w \in C. \tag{3.5}$$

From condition a1: $AC = D$ and by (3.5), we have

$$\langle g(Ax), Au\rangle \leq \langle g(Ax), z\rangle, \text{ for all } z \in D. \tag{3.6}$$

It implies that

$$Au \in G(Ax). \tag{3.7}$$

Hence $(u, Au) \in T(x, Ax)$; and therefore, $T(x, Ax) \neq \emptyset$.

To prove that $T: M \to 2^M \setminus \{\emptyset\}$ is $\succcurlyeq^{X \times Y}$-increasing upward, we take arbitrary $(t, At), (s, As) \in M$ with $(s, As) \succcurlyeq^{X \times Y} (t, At)$ that is equivalent to $s \succcurlyeq^X t$. It is because that $A$ is order-increasing. For any given $(p, Ap) \in T(t, At)$, we have $p \in F(t)$ and $Ap \in G(At)$. From Lemma 3.6, and the chain complete property of $(F(t), \succcurlyeq^X)$, we can choose an $\succcurlyeq^X$-maximal point $p' \in F(t)$ such that $p' \succcurlyeq^X p$ (From Lemma 3.6, $p'$ is also an $\succcurlyeq^X$-maximal point in $(C, \succcurlyeq^X)$). From Lemma 3.5, $F$ is $\succcurlyeq^X$-increasing upward. By $s \succcurlyeq^X t$ and $p' \in F(t)$, there is $q \in F(s)$ such that $q \succcurlyeq^X p'$. We can similarly choose an $\succcurlyeq^X$-maximal point $q'$ in $(F(s), \succcurlyeq^X)$ such that $q' \succcurlyeq^X q$. Then, similarly to the proofs of (3.5), (3.6) and (3.7), we can show that

$$Aq' \in G(As) \text{ (and } Ap' \in G(As)).$$

Then we obtain

$$(q', Aq') \in T(s, As). \tag{3.8}$$

From $q' \succcurlyeq^X q \succcurlyeq^X p' \succcurlyeq^X p$, and by the order-increasing property of $A$, it implies $Aq' \succcurlyeq^Y Aq \succcurlyeq^Y Ap' \succcurlyeq^Y Ap$. It follows that

$$(q', Aq') \succcurlyeq^{X \times Y} (p, Ap). \tag{3.9}$$

By (3.8) and (3.9), it implies that $T: M \to 2^M \setminus \{\emptyset\}$ is $\succcurlyeq^{X \times Y}$-increasing upward.

Since $T(x, Ax)$ is a compact subset in $C \times D$ with respect to the product topology, by Lemma 5.3 [13], $(T(x, Ax), \succcurlyeq^{X \times Y})$ is universally inductive, for every $(x, Ax) \in M$.

From the order increasing property of $A$, it can be checked that the points $(x_0, Ax_0) \in M$ and $(w, Aw) \in F(x_0) \times G(Ax_0)$ with $x_0 \preccurlyeq^X w$, given in condition (3.2) in this theorem satisfy

$$(w, Aw) \in T(x_0, Ax_0) \text{ with } (x_0, Ax_0) \preccurlyeq^{X \times Y} (w, Aw).$$

Hence, $T$ satisfies all conditions A1, A2, and A3 in Theorem 3.1 [13]. Then $T$ has a fixed point $(x^*, Ax^*)$ in $M$, such that $(x^*, Ax^*) \in T(x^*, Ax^*)$. It implies

$$x^* \in F(x^*) = \{u \in C: \langle f(x^*), u\rangle = \min_{t \in C} \langle f(x^*), t\rangle\},$$

and

$$Ax^* \in G(Ax^*) = \{v \in D: \langle g(Ax^*), v\rangle = \min_{s \in C} \langle g(Ax^*), s\rangle\}.$$

It follows that

$\langle f(x^*), x - x^* \rangle \geq 0$, for all $x \in C$ and $\langle g(Ax^*), y - Ax^* \rangle \geq 0$, for all $y \in D$.

Hence $(x^*, Ax^*)$ is a solution of the split variational inequality problem SVIP($f, C, A, g, D$). The conclusions (i) and (ii) in this theorem immediately follow from the conclusions (i) and (ii) in Theorem 3.1 [13]. □

When the considered partially ordered Banach spaces in split variational inequality problems are reflexive, we have the following result.

**Theorem 3. 9**. *Let $(X, \succcurlyeq^X)$ and $(Y, \succcurlyeq^Y)$ be partially ordered reflexive Banach spaces and let $C \subseteq X$ and $D \subseteq Y$ be nonempty, bounded, closed and convex subsets. Let $f: C \to X^*$, $g: D \to Y^*$ be continuous operators and let $A: X \to Y$ be a bounded linear operator satisfying the following conditions*:

    (a0) *$f$ is $\succcurlyeq^X$-decrement and $g$ is $\succcurlyeq^Y$-decrement*;
    (a1) *$AC = D$*;
    (b1) *$A$ is order-increasing*;
    (b2) *$A$ is order-pseudomonotone with respect to $f$ and $g$*.

*Suppose that there is $(x_0, Ax_0) \in C \times D$ such that*

$$\text{there exists } (w, Aw) \in F(x_0) \times G(Ax_0) \text{ satisfying } x_0 \preccurlyeq^X w. \tag{3.2}$$

*Then* SVIP($f, C, A, g, D$) *has a solution. Moreover, we have*

    (i)   $(\mathcal{S}(f, C, A, g, D), \succcurlyeq^{X \times Y})$ *is a nonempty inductive poset*;
    (ii) $(\mathcal{S}(f, C, A, g, D) \cap [(x_0, Ax_0)), \succcurlyeq^{X \times Y})$ *is also a nonempty inductive poset*.

**Sketch of the proof**. Since $C$ and $D$ are nonempty, bounded, closed and convex subsets of partially ordered reflexive Banach spaces $(X, \succcurlyeq^X)$ and $(Y, \succcurlyeq^Y)$, then $C$ and $D$ are weakly compact subsets in $X$ and $Y$, respectively. The graph $M$ of the continuous operator $A$ is a weakly compact subset in $C \times D$ with respect to the product weak topology. From Observation 3.1, and applying Lemma 3.2, we obtain that $(M, \succcurlyeq^{X \times Y})$ is an $\succcurlyeq^{X \times Y}$-chain complete poset.

Since $f$ and $g$ are continuous, from Lemma 3.3, the minimizing acceptances $F(x)$ and $G(y)$ are closed subsets in the weakly compact sets $C$ and $D$, respectively, which implies that $F(x)$ and $G(y)$ are weakly compact subsets in $X$ and $Y$, for every $x \in C$ and $y \in D$, respectively. Hence, both $(F(x), \succcurlyeq^X)$ and $(G(y), \succcurlyeq^Y)$ are chain complete. It implies that the operator $T$ defined by (3.3) has $\succcurlyeq^{X \times Y}$-chain complete values, that is, $(T(x, Ax), \succcurlyeq^{X \times Y})$ is an $\succcurlyeq^{X \times Y}$-chain complete subset in the $\succcurlyeq^{X \times Y}$-chain complete poset $(M, \succcurlyeq^{X \times Y})$, for all $(x, Ax) \in M$.

Then rest of the proof is similar to the proof of Theorem 3.8, by applying Theorem 3.1 [13]. □

As a consequence of Theorems 3.8 and 3.9, from the properties of inductive posets, we can obtain the existence of order maximal solutions to some split variational inequality problems.

**Remarks 3.10**. Under the hypotheses of Theorem 3.8 or Theorem 3.8, the split variational inequality problem has the following property:

(I) SVIP($f$, $C$, $A$, $g$, $D$) has an $\succcurlyeq^{X \times Y}$-maximal solution;
(II) SVIP($f$, $C$, $A$, $g$, $D$) has an $\succcurlyeq^{X \times Y}$-maximal solution $(x^*, Ax^*)$ such that $x^* \succcurlyeq^X x_0$ and $Ax^* \succcurlyeq^Y Ax_0$, where the point $(x_0, Ax_0)$ is given in (3.2).

### 4. Split convex minimization problems

Let $X$ and $Y$ be Banach spaces and let $(U, \succcurlyeq)$ be a partially ordered vector space. For given nonempty, closed and convex subsets $C \subseteq X$ and $D \subseteq Y$, given operators $\varphi: C \to U$, $\psi: D \to U$ and a bounded linear operator $A: X \to Y$, the ordered split convex minimization problem associated with $\varphi$, $C$, $A$, $\psi$, $D$, $U$ is formulated as follows:

$$\text{find } x_* \in C \text{ such that } \varphi(x_*) \preccurlyeq \varphi(x), \text{ for all } x \in C, \tag{4.1}$$

$$\text{and such that the point } y_* = Ax_* \in D \text{ solves } \psi(y_*) \preccurlyeq \psi(y), \text{ for all } y \in D. \tag{4.2}$$

Such a pair $(x_*, y_*) = (x_*, Ax_*)$ is called a solution of the above ordered split convex minimization problem.

An operator $\varphi: C \to U$ is said to be $\succcurlyeq$-continuous if, for any $u \in U$, the following subsets are closed

$$\{x \in C: \varphi(x) \preccurlyeq u\} \text{ and } \{x \in C: \varphi(x) \succcurlyeq u\}.$$

As a special case of ordered split convex minimization problems, if $U$ is taken to be the real set $R$ with the ordinary order, then the ordered split convex minimization problem defined in (4.1) and (4.2) becomes the ordinary split convex minimization problem:

$$\text{find } x_* \in C \text{ such that } \varphi(x_*) \leq \varphi(x), \text{ for all } x \in C, \tag{4.3}$$

$$\text{and such that the point } y_* = Ax_* \in D \text{ solves } \psi(y_*) \leq \psi(y), \text{ for all } y \in D. \tag{4.4}$$

A bounded linear operator $A: X \to Y$ is said to have the *convexity order reserving property* (or it is said to be *convexity order reserved*) with respect to the mappings $\varphi$ and $\psi$ on $C$ and $D$ if, for any given vectors $x_1, x_2, \cdots, x_n$ in $C$, and for any strictly linear combination $x' = \Sigma_{1 \leq i \leq n} \lambda_i x_i$, there is a positive integer $j \leq n$ such that

$$\varphi(x') \preccurlyeq \varphi(x_j) \text{ and } \psi(Ax') \preccurlyeq \psi(Ax_j).$$

In particular, if the partially ordered vector space $(U, \succcurlyeq)$ is taken to be $(R, \geq)$, then the bounded linear operator $A$ is said to have *convexity value reserving property* (or it is said to be *convexity value reserved*) with respect to the mappings $\varphi$ and $\psi$ on $C$ and $D$.

**Theorem 4.1**. *Let $X$ and $Y$ be Banach spaces and $C \subseteq X$ and $D \subseteq Y$ be nonempty, closed and convex subsets. Let $(U, \succcurlyeq)$ be a partially ordered vector space. Given $\succcurlyeq$-continuous operators $\varphi: C \to U$, $\psi: D \to U$ and a bounded linear operator $A: X \to Y$ satisfying the following conditions*:

       (a1) $AC = D$;
       (a2)' $A$ is convexity order reserved with respect to $\varphi$ and $\psi$ on $C$ and $D$.

*Suppose that there is a point $(x_0, Ax_0) \in C \times D$ such that the set*

$$\{(x', Ax') \in C \times D: \varphi(x') \preccurlyeq \varphi(x_0) \text{ and } \psi(Ax') \preccurlyeq \psi(Ax_0)\} \text{ is compact.}$$

*Then the ordered split convex minimization problem defined by (4.1) and (4.2) has a solution.*

**Proof**. For the given partially ordered vector space $(U, \succcurlyeq)$, we define the component partial order $\succcurlyeq^2$ on the product space $U \times U$ as follows: for any $(u_1, v_1), (u_2, v_2) \in U \times U$,

$$(u_2, v_2) \succcurlyeq^2 (u_1, v_1) \text{ if and only if } u_2 \succcurlyeq u_1 \text{ and } v_2 \succcurlyeq v_1.$$

It follows that $(U \times U, \succcurlyeq^2)$ is also a partially ordered vector space. Based on the mappings $\varphi$ and $\psi$, we define a (product) mapping $\mu: C \times D \to U \times U$ as

$$\mu(x, y) = (\varphi(x), \psi(y)), \text{ for every } (x, y) \in C \times D. \tag{4.5}$$

Then, by using the mapping $\mu$ in (4.5), the ordered split minimization problem defined by (4.1) and (4.2) can be transformed into the following ordered convex minimization problem:

find a point $(x^*, Ax^*) \in C \times D$ such that $\mu(x^*, Ax^*) \preccurlyeq^2 \mu(x, y)$, for all $(x, y) \in C \times D$.

The remainder of the proof is similar to the proof of Theorem 2.3. □

As an immediate consequence of Theorem 4.1, if we take the partially ordered vector space $U$ in Theorem 4.1 to be the real numbers set $R$, we have the existence of solutions to the split convex minimization problem defined by (4.3) and (4.4).

**Corollary 4.2**. *Let $X$ and $Y$ be Banach spaces and $C \subseteq X$ and $D \subseteq Y$ be nonempty, closed and convex subsets. Given continuous functions $\varphi: C \to R$, $\psi: D \to R$ and a bounded linear operator $A: X \to Y$ satisfying the following conditions*:

       (a1) $AC = D$;
       (a2)'' $A$ is convexity value reserved with respect to $\varphi$ and $\psi$ on $C$ and $D$.

*Suppose that there is a point* $(x_0, Ax_0) \in C \times D$ *such that*

$$\{(x', Ax') \in C \times D: \varphi(x') \leq \varphi(x_0) \text{ and } \psi(Ax') \leq \psi(Ax_0)\} \text{ is compact}.$$

*Then the split convex minimization problem defined by* (4.3) *and* (4.4) *has a solution*.

**Acknowledgements:** The author is very grateful to Professor Hong-kun Xu for his kind encouragements and valuable suggestions for the research of this paper.